\documentstyle[11pt,amssymb]{article}

\newtheorem{prop}{Proposition}
\newtheorem{thm}{Theorem}
\newtheorem{cor}{Corollary}

\def\Z{{\Bbb Z}}
\def\Q{{\Bbb Q}}
\def\R{{\Bbb R}}
\def\C{{\Bbb C}}
\def\O{{\cal O}}
\def\E{{\cal E}}
\def\H{{\cal H}}
\def\L{{\Lambda}}
\def\om{{\omega}}
\def\endproof{\hfill $\Box$}

\newcommand{\dis}{\displaystyle}

\newcommand{\hin}{\hspace{-2pt}\in\hspace{-2pt}}

\newcommand{\gequ}{\geqslant}
\newcommand{\lequ}{\leqslant}
\newcommand{\ra}{\rightarrow}

\newcommand{\dbar}{\overline{\partial}}
\newcommand{\Spec}{\mbox{Spec}}

\newcommand{\Hom}{\mbox{Hom}}
\newcommand{\End}{\mbox{End}}

\newcommand{\Ker}{\mbox{Ker}}
\newcommand{\Tr}{\mbox{Tr}}

\newcommand{\ov}{\overline}
\newcommand{\noin}{\noindent}
\newcommand{\wt}{\widetilde}
\newcommand{\wh}{\widehat}
\newcommand{\pr}{\prime}

\newcommand{\rk}{\mbox{rk}}

\begin{document}

\title{\bf Bott-Chern Forms and Arithmetic Intersections}
\author{Harry Tamvakis\\
Department of Mathematics\\
The University of Chicago\\
Chicago, IL 60637}
\date{}
\maketitle 

\begin{abstract}
Let  $\ov{\E}:\ 0\ra \ov{S}\ra\ov{E}\ra\ov{Q}\ra 0$ be
a short exact sequence of hermitian
vector bundles with metrics on $S$ and $Q$ induced from that on $E$.
We compute the Bott-Chern form  $\wt{\phi}(\ov{\E})$ corresponding to
any characteristic class $\phi$, assuming $\ov{E}$ is projectively flat. The
result is used to obtain a new presentation of the Arakelov Chow ring of the
arithmetic grassmannian.
\end{abstract}

\section{Introduction}

  Arakelov theory is an intersection theory for varieties
over rings $\O_F$ of algebraic integers, analogous to the usual
one over fields. The fundamental idea is that in order to have a good 
theory of intersection numbers, one has to include
information at the infinite primes.
  
  The work of Arakelov in dimension two has been generalized by
Gillet and Soul\'{e} to higher dimensional {\em arithmetic varieties} $X$,
by which we mean  regular, projective and flat schemes over $\Spec\Z$. They
define an {\em arithmetic Chow ring} $\wh{CH}(X)_{\Q}$ whose elements are 
represented
by cycles on $X$ together with Green currents on $X(\C)$.
The theory is a blend of arithmetic, algebraic
geometry and complex hermitian geometry. For example, the
Faltings height of an arithmetic variety $X$ is realized as an
`arithmetic degree' with respect to a hermitian line bundle
over $X$. 

  A {\em hermitian vector bundle} $\ov{E}=(E,h)$ over $X$ is an algebraic
vector bundle $E$ on $X$ together with a hermitian metric $h$ on the
corresponding holomorphic vector bundle $E(\C)$ on the complex
manifold $X(\C)$. To such an object one associates arithmetic
Chern classes $\wh{c}(\ov{E})$
with values in $\wh{CH}(X)$. These satisfy most of the
usual properties of Chern classes, with one exception:
 given a short exact sequence of hermitian vector bundles
\begin{equation}
\label{seq}
\ov{\E}:\ 0\ra \ov{S}\ra\ov{E}\ra\ov{Q}\ra 0
\end{equation}
the class $\wh{c}(\ov{S})\wh{c}(\ov{Q})-\wh{c}(\ov{E})$ vanishes
 when $\ov{E}$ is the orthogonal direct sum of $\ov{S}$ and
$\ov{Q}$. In general however this difference is non-zero and 
is the image in
$\wh{CH}(X)$ of a differential form on $X(\C)$, the {\em Bott-Chern
form} associated to the exact sequence $\ov{\E}$.

 These secondary characteristic classes
were originally defined by Bott and Chern [BC]
with applications to value distribution theory. They later occured
in the work of Donaldson [Do] on Hermitian-Einstein metrics.
Bismut, Gillet and Soul\'{e} [BiGS] gave a new axiomatic definition for 
Bott-Chern forms, suitable for use in arithmetic intersection theory.
Given an exact sequence of hermitian 
holomorphic vector bundles as in (\ref{seq}),
we have
$c(\ov{S})c(\ov{Q})-c(\ov{E})=dd^c\eta$ for some form $\eta$; the
Bott Chern form of $\ov{\E}$ is a natural choice of such an $\eta$.

  Calculating these forms is important because they give relations
in the arithmetic Chow ring of an arithmetic variety. No systematic
work has appeared on this; rather one finds scattered
calculations throughout the literature (see for example 
[BC], [C1], [D], [GS2], [GSZ], [Ma], [Mo]).
We confine ourselves to the case where
the metrics on $S$ and $Q$ are induced from the one on $E$.
Our goal is to give explicit
formulas for the Bott-Chern forms corresponding to {\em any} characteristic
class, when they can be expressed
in terms of the characteristic classes of the bundles involved. 
This is not always possible as these forms are not closed in
general; however the situation is completely understood 
when $E$ is a projectively flat bundle. The results
  build on the work of Bott, Chern, Cowen, Deligne, Gillet,
Soul\'{e} and Maillot. Some of our calculations overlap with previous work,
but with simpler proofs. 

The main application we give to arithmetic
intersection theory is a new presentation of the
 Arakelov Chow ring of the grassmannian over $\Spec \Z$. Maillot [Ma] gave a
presentation of this ring and formulated an `arithmetic Schubert calculus'.
We hope our work contributes towards a better understanding of these 
intersections.

This paper is organized as follows. Section \ref{isf} is a review
of some basic material on invariant and symmetric functions.
In \S \ref{bcfs} we recall the hermitian geometry we will need,
including the definition of Bott-Chern forms.
The basic tool for calculating these forms is reviewed in
\S \ref{cbcf}, with some applications that have appeared
before in the literature. 
\S \ref{ait} is mainly 
an exposition of the arithmetic intersection theory that we require. 
The rest of the paper is new.
In sections \ref{flatses} and \ref{projflat}
we derive formulas for computing Bott-Chern forms of short exact sequences
(with the induced metrics)
for any characteristic class when $\ov{E}$ is flat or more generally
projectively flat.
  We emphasize the central role played by the
{\em power sum forms} in the results; to our knowledge this
phenomenon has not been observed before.
 The combinatorial identities involving harmonic numbers that we 
encounter are also interesting. Sections 2-6 contain results in hermitian
complex geometry and may be read without prior knowledge of Arakelov theory.
\S \ref{grass} applies our calculations to obtain a 
presentation of the Arakelov Chow ring of the arithmetic grassmannian.

This should be regarded as a companion paper to [T]; both papers will be
 part of
the author's 1997 University of Chicago thesis. I wish to thank my 
advisor William Fulton for many useful conversations and 
exchanges of ideas.

\section{Invariant and symmetric functions}
\label{isf}

The symmetric group $S_n$ acts on the polynomial ring
$\Z[x_1,x_2,\ldots,x_n]$ by permuting the variables, and the ring of 
invariants $\L(n)=\Z[x_1,x_2,\ldots,x_n]^{S_n}$ is the ring of symmetric
polynomials. For $B=\Q$ or $\C$, let $\L(n,B)=\L(n)\otimes_{\Z} B$. 

Let $e_k(x_1,\ldots,x_n)$ be the $k$-th elementary symmetric polynomial
in the variables $x_1,\ldots,x_n$
and $\dis p_k(x_1,\ldots,x_n)=\sum_i x_i^k$ the $k$-th power sum. 
The fundamental theorem on symmetric functions states that 
$\L(n)=\Z[e_1,\ldots,e_n]$ and that $e_1,\ldots,e_n$ are algebraically
independent.
For $\lambda$ a partition, i.e. a decreasing sequence 
$\lambda_1\gequ \lambda_2\gequ\cdots\gequ\lambda_m$ of nonnegative
integers, define $\dis p_{\lambda}:=\prod_{i=1}^mp_{\lambda_i}$. 
It is well known that the $p_{\lambda}$'s form an additive $\Q$-basis 
for the ring of symmetric polynomials (cf. [M], \S 2). The two bases
are related by Newton's identity:
\begin{equation}
\label{ni}
p_k-e_1p_{k-1}+e_2p_{k-2}-\cdots+(-1)^kke_k=0.
\end{equation}

Another important set of symmetric functions related to the cohomology
ring of grassmannians are the Schur polynomials. For a partition $\lambda$
as above, the Schur polynomial $s_{\lambda}$ is defined by
\[
\dis
s_{\lambda}(x_1,\ldots,x_n)=\frac{1}{\Delta}\cdot
\det(x_i^{\lambda_j+n-j})_{1\lequ i,j\lequ n},
\]
where $\dis\Delta=\prod_{1\lequ i<j\lequ n}(x_i-x_j)$ is the Vandermonde
determinant. The $s_{\lambda}$ for all $\lambda$ of length $m\lequ n$ form
a $\Z$-basis of $\L(n)$ (cf. [M], \S I.3). 

Let $\C[T_{ij}]$ $(1\lequ i,j\lequ n$) be the coordinate ring of the 
space $M_n(\C)$ of $n\times n$ matrices. $GL_n(\C)$ acts on matrices by 
conjugation, and we let $I(n)=\C[T_{ij}]^{GL_n(\C)}$ denote the 
corresponding graded ring of invariants. There is an isomorphism
$\tau : I(n)\ra \L(n,\C)$ given by evaluating
an invariant polynomial $\phi$ on the diagonal matrix 
diag$(x_1,\ldots,x_n)$. We will often identify $\phi$ with the the
symmetric polynomial $\tau(\phi)$. 
We will need to consider invariant
polynomials with rational coefficients; let
$I(n,\Q)\simeq\Q[x_1,x_2,\ldots,x_n]^{S_n}$ be the corresponding ring.

Given $\phi\hin I(n)$, let $\phi^{\pr}$ be a $k$-multilinear form on $M_n(\C)$
 such that
\[
\dis
\phi^{\pr}(gA_1g^{-1},\ldots,gA_kg^{-1})=\phi^{\pr}(A_1,\ldots,A_k)
\]
for $g\hin GL(n,\C)$ and $\phi(A)=\phi^{\pr}(A,A,\ldots,A)$. 
Such forms are most easily constructed for the power sums $p_k$ by setting
\[
p_k^{\pr}(A_1,A_2,\ldots,A_k)=\Tr (A_1A_2\cdots A_k).
\]
For
$p_{\lambda}$ we can take $p_{\lambda}^{\pr}=\prod p_{\lambda_i}^{\pr}$.
Since the $p_{\lambda}$'s are a basis 
of $\L(n,\Q)$, it follows that
 one can use the above constructions to find 
multilinear forms $\phi^{\pr}$ for any $\phi\hin I(n)$.

An explicit formula for $\phi^{\pr}$ is given by polarizing $\phi$:
\[
\dis
\phi^{\pr}(A_1\ldots,A_k)=\frac{(-1)^k}{k!}\sum_{j=1}^k
\sum_{i_1<\ldots<i_j}(-1)^j\phi(A_{i_1}+\ldots+A_{i_j}).
\]
Although above formula for $\phi^{\pr}$ is symmetric in $A_1\ldots,A_k$, this
property is not needed for the applications that follow.

\section{Hermitian differential geometry}
\label{bcfs}

 Let $X$ be a complex manifold, $E$ a rank $n$ holomorphic vector
bundle over $X$. Denote by $A^k(X,E)$ the $C^{\infty}$ sections of 
$\Lambda^kT^*X\otimes E$, where $T^*X$ denotes the cotangent bundle of 
$X$. In particular $A^k(X)$ is the space of smooth complex $k$-forms on $X$.
Let $A^{p,q}(X)$ the space of smooth complex forms of type $(p,q)$ on
$X$ and $A(X):=\bigoplus_pA^{p,p}(X)$.
 The decomposition $A^1(X,E)=A^{1,0}(X,E)\bigoplus A^{0,1}(X,E)$
induces a decomposition $D=D^{1,0}+D^{0,1}$ of each connection $D$ on $E$.
 Let $d=\partial+\overline{\partial}$
and $d^c=(\partial-\overline{\partial})/(4\pi i)$.

 Assume now that $E$ is equipped with a hermitian metric $h$. The pair 
$(E,h)$ is called a {\em hermitian vector bundle}. The metric $h$ induces a 
canonical connection $D=D(h)$ such that $D^{0,1}=\overline{\partial}_E$
and $D$ is {\em unitary}, i.e. 
\[
\dis
d\,h(s,t)=h(Ds,t)+h(s,Dt), 
\mbox{ for all } s,t\hin A^0(X,E).
\]
 The connection $D$ is called the {\em hermitian
holomorphic connection} of $(E,h)$. $D$ can be extended to $E$-valued forms
by using the Leibnitz rule:
\[
\dis
D(\om\otimes s)=d\om\otimes s+ (-1)^{\deg\om}\om\otimes Ds.
\]
 The composite
\[
\dis
K=D^2:A^0(X,E)\ra A^2(X,E)
\]
is $A^0(X)$-linear; hence $K\hin A^2(X,\End(E))$. In fact
$K=D^{1,1}\hin
A^{1,1}(X,\End(E))$, because $D^{0,2}=\dbar^2_E=0$, so $D^{2,0}$ also
vanishes by unitarity. $K$ is called the {\em curvature} of $D$.

Given a hermitian vector bundle $\ov{E}=(E,h)$ and an invariant polynomial
 $\phi\hin I(n)$ there is an associated differential form 
$\phi(\ov{E}):=\phi(\frac{i}{2\pi}K)$,
 defined locally by identifying $\End(E)$ with $M_n(\C)$; 
$\phi(\ov{E})$ makes sense globally on $X$ since $\phi$ is invariant by
conjugation. These differential forms are $d$ and $d^c$ closed and
have the following properties (cf. [BC]):

\noindent
(i) The de Rham cohomology class of $\phi(\ov{E})$ is independent of the 
metric $h$ and coincides with the usual characteristic class from topology.

\noindent
(ii) For every holomorphic map $f:X\ra Y$ of complex manifolds,
\[
f^*(\phi(E,h))=\phi(f^*E,f^*h).
\]

One thus obtains the {\em Chern forms} $c_k(\ov{E})$ with $c_k=e_k(x_1,
\ldots,x_n)$, the {\em power sum forms} $p_k(\ov{E})$, the
{\em Chern character form} $ch(\ov{E})$ with 
$\dis ch(x_1,\ldots,x_n)=\sum_i \exp(x_i)=\sum_k \frac{1}{k!}p_k$, etc.

We fix some more notation: A direct sum $\ov{E}_1\bigoplus\ov{E}_2$ of 
hermitian vector bundles will always
mean the orthogonal direct sum
$(E_1\bigoplus E_2,h_1\oplus h_2)$. 
Let $\wt{A}(X)$ be the
quotient of $A(X)$ by $\mbox{Im} \partial + \mbox{Im} \dbar$. 
If $\omega$ is a closed form in $A(X)$ the cup product $\wedge\omega:
\wt{A}(X)\ra\wt{A}(X)$ and the operator
$dd^c:\wt{A}(X)\ra A(X)$ are well defined.

 Let $\E:\ 0\ra S \ra E\ra Q \ra 0$ be an exact 
sequence of holomorphic vector bundles on $X$. Choose arbitrary hermitian
metrics $h_S,h_E,h_Q$ on $S,E,Q$ respectively.
Let
\[
\dis
\ov{\E}=(\E,h_S,h_E,h_Q):\
0\ra \ov{S} \ra \ov{E}\ra \ov{Q}\ra 0.
\]
Note that we do not in general assume that the metrics $h_S$ or
$h_Q$ are induced from $h_E$. We say that $\ov{\E}$ is {\em split}
when $(E,h_E)=(S\bigoplus Q,h_S\oplus h_Q)$ and 
$\E$ is the obvious exact sequence. Following [GS2], we have the
following

\begin{thm} \label{bc}
 Let $\phi\hin I(n)$ be any invariant polynomial.
 There is a unique way to attach to every exact
sequence $\ov{\E}$ a form $\wt{\phi}(\ov{\E})$ in 
$\wt{A}(X)$ in such a way that:

\noin
{\em (i)} $dd^c\wt{\phi}(\ov{\E})=\phi(\ov{S}\bigoplus \ov{Q})
-\phi(\ov{E})$,

\noin
{\em (ii)} For every map $f:X\ra Y$ of complex manifolds,
        $\wt{\phi}(f^*(\ov{\E}))=f^*\wt{\phi}(\ov{\E})$,

\noin
{\em (iii)} If $\ov{\E}$ is split, then $\wt{\phi}(\ov{\E})=0$.

\end{thm}

 In [BC], Bott and Chern solved the equation 
$dd^c\wt{\phi}(\ov{\E})=\phi(\ov{S}\bigoplus \ov{Q})
-\phi(\ov{E})$ when the metrics on $S$ and $Q$ are induced from the
metric on $E$. In [BiGS] a new axiomatic definition of these forms was given,
more generally for an acyclic complex of holomorphic vector bundles on $X$.

The following useful calculation is an immediate consequence of the
definition ([GS2], Prop. 1.3.1):

\begin{prop} \label{bcprop}
 Let $\phi$ and $\psi$ be two invariant polynomials. Then
\[
\dis
\wt{\phi + \psi}(\ov{\E})=\wt{\phi}(\ov{\E})+\wt{\psi}(\ov{\E}).
\]
\[
\wt{\phi\psi}(\ov{\E})=\wt{\phi}(\ov{\E})\psi(\ov{E})+
\phi(\ov{S}\oplus \ov{Q})\wt{\psi}(\ov{\E})=
\wt{\phi}(\ov{\E})\psi(\ov{S}\oplus \ov{Q})+
\phi(\ov{E})\wt{\psi}(\ov{\E}).
\]

\end{prop}
{\bf Proof.} One checks that right hand side of these identities 
satisfies the
three properties of Theorem \ref{bc} that characterize the left hand side.
\endproof

\medskip

We will also need to know the behaviour of $\wt{c}$ when $\ov{\E}$ is
twisted by a line bundle. The following is a consequence of
[GS2], Prop. 1.3.3:

\begin{prop}
\label{twist}
For any hermitian line bundle $\ov{L}$,
\[
\dis
\wt{c_k}(\ov{\E}\otimes\ov{L})=
\sum_{i=1}^k{n-i \choose k-i}\wt{c_i}(\ov{\E})c_1(\ov{L})^{k-i}.
\]
\end{prop}

\section{Calculating Bott-Chern Forms}
\label{cbcf}

In this section we will consider an exact sequence 
\[
\dis
\ov{\E}:\
0\ra \ov{S} \ra \ov{E}\ra \ov{Q}\ra 0.
\]
where the metrics on $\ov{S}$ and $\ov{Q}$ are induced from the
metric on $E$. Let $r$, $n$ be the ranks of the bundles $S$ and $E$.
Let $\phi\hin I(n)$ be homogeneous of degree $k$.
We will formulate a theorem for calculating the Bott-Chern form
$\wt{\phi}(\ov{\E})$. This result follows from the work of
Bott-Chern, M. Cowen, J. Bismut and Gillet-Soul\'{e}.

Let $\phi^{\pr}$ be defined as in \S \ref{isf}.
For any two matrices $A,B \hin M_n(\C)$ set
\[
\dis
\phi^{\pr}(A;B):=\sum_{i=1}^k\phi^{\pr}(A,A,\ldots,A,B_{(i)},A,\ldots,A),
\]
where the index $i$ means that $B$ is in the $i$-th position.

Choose a local orthonormal frame $s=(s_1,s_2,\ldots,s_n)$
of $E$ such that the first $r$ elements generate $S$, and
let $K(\ov{S})$, $K(\ov{E})$ and $K(\ov{Q})$ 
be the curvature matrices of $\ov{S}$, $\ov{Q}$
and $\ov{E}$ with respect to $s$.
Let $K_S=\frac{i}{2\pi}K(\ov{S})$,
$K_E=\frac{i}{2\pi}K(\ov{E})$ and
$K_Q=\frac{i}{2\pi}K(\ov{Q})$. The matrix $K_E$ has the form
\[
\dis
K_E=
\left(
\begin{array}{c|c}
K_{11} & K_{12} \\ \hline
K_{21} & K_{22} 
\end{array} \right)
\]
where $K_{11}$ is an $r\times r$ submatrix. Also consider the matrices
\[
\dis
K_0=
\left(
\begin{array}{c|c}
K_S & 0 \\ \hline
K_{21} & K_Q 
\end{array} \right)
\ \mbox{ and } \
J_r=
\left(
\begin{array}{c|c}
Id_r & 0 \\ \hline
0 & 0 
\end{array} \right).
\]

 Let $u$ be a formal variable and
$K(u):=uK_E+(1-u)K_0$.
Finally, let $\phi^!(u)=\phi^{\pr}(K(u); J_r)$.
 We then have the following

\begin{thm} \label{calc}

\begin{equation}
\label{cowen}
\wt{\phi}(\ov{\E})=\int_0^1\frac{\phi^!(u)-\phi^!(0)}{u}\,du.
\end{equation}

\end{thm}

\medskip
{\bf Proof.} We prove that $\wt{\phi}(\ov{\E})$ as defined above satisfies
axioms (i)-(iii) of Theorem \ref{bc}. The main step is the first
axiom; this was essentially done in [BC]\, \S 4, when $\phi=c$ is
the total Chern class. In the form (\ref{cowen}) 
(again for the total Chern class), the equation was given
by M. J. Cowen in [C1] and [C2], while simplifying Bott and Chern's
proof. We follow both sources in sketching a proof of this more general result.

Let $h$ and $h_Q$ denote the metrics on $E$ and $Q$ respectively.
Define the orthogonal projections $P_1:\ov{E}\ra \ov{S}$ and 
$P_2:\ov{E}\ra \ov{Q}$ and put 
$h_u(v,v^{\pr})=uh(P_1v,P_1v^{\pr})+h(P_2v,P_2v^{\pr})$ for 
$v,v^{\pr}\hin E_x$ and $0 < u\lequ 1$.
Then $h_u$ is a hermitian norm, $h_1=h$ and $h_u\ra h_Q$ as $u\ra 0$.
Let $K(E,h_u)$ be the curvature matrix of $(E,h_u)$
relative to the holomorphic frame $s$ defined above. 
 Proposition 3.1 of [C2] proves that $\frac{i}{2\pi}K(E,h_u)=K(u)$.
It follows from Proposition 3.28 of [BC] that for $0 < t \lequ 1$,
\[
\dis
\phi(E,h_t)-\phi(E,h)=
dd^c\int_t^1
\frac{\phi^{\pr}(K(u); J_r)}{u}\,du.
\]

If we could set $t=0$ we would be done; however, the integral will not be
convergent in general. Note that $K(u)=K_0+uK_1$, where $K_1\hin
A^{1,1}(X,\End(E))$ is independent of $u$. Therefore it will suffice
to show that $\phi^{\pr}(K_0; J_r)$ is a closed form, so that it can be
deleted from the integral. For this we may assume that $\phi=p_{\lambda}$ is
a product of power sums, $\lambda=(\lambda_1,\lambda_2,\ldots,\lambda_m)$
a partition. Then
\[
\dis
p_{\lambda}^{\pr}(K_0; J_r)=
\sum_{i=1}^m  \Tr(K_S)^{\lambda_i-1} \prod_{j\neq i}
(\Tr(K_S)+\Tr(K_Q))^{\lambda_j}=
\sum_{i=1}^m p_{\lambda_i-1}(\ov{S})\prod_{j\neq i}
p_{\lambda_j}(\ov{S}\oplus\ov{Q})
\]
is certainly a closed form.

This proves axioms (i) and (iii); axiom (ii) is easily checked as well.
\endproof

\medskip

\noindent
{\bf Remark.}
A similar deformation to the one in [C2] was used 
by Deligne in
[D], 5.11 for a calculation involving the Chern character form.
Special cases of Theorem \ref{calc} have been used in the literature before, 
see for example [GS2] Prop. 5.3, [GSZ] 2.2.3 and [Ma] Theorem 3.3.1.

\medskip

We deduce some simple but useful calculations:

\begin{cor} \label{c1}
{\em (a)} $\wt{c_1^k}(\ov{\E})=0$ for all $k\gequ 1$ and 
$\wt{c}_m(\ov{\E})=0$ for all $m > ${\em rk}$ E$.

\noindent
{\em (b)} $\wt{p_2}(\ov{\E})=2(${\em Tr}$K_{11}-c_1(\ov{S}))$ and
$\wt{c_2}(\ov{\E})=c_1(\ov{S})-${\em Tr}$K_{11}$.
 
\end{cor}
{\bf Proof.} (a) $c^!_1(u)$ is
independent of $u$; hence $\wt{c_1}(\ov{\E})=0$. The result for
higher powers of $c_1$ follows from Proposition 
\ref{bcprop}. In addition, $\wt{c}_m(\ov{\E})=0$ for $m > \rk E$ is an
immediate consequence of the definition.

\noindent
(b) Using the bilinear form $p_2^{\pr}$ described previously,
we find $p^!_2(u)=2(u\Tr K_{11}+(1-u)c_1(\ov{S}))$, so
\[
\dis
\wt{p_2}(\ov{\E})=2\int_0^1\frac{u\Tr K_{11}
+(1-u)c_1(\ov{S})-c_1(\ov{S})}{u}\,du=2(\Tr K_{11}-c_1(\ov{S})).
\]
To calculate $\wt{c_2}(\ov{\E})$, use the identity 
$2c_2=c_1^2-p_2$.\
\endproof

\medskip
 
 Corollary \ref{c1}(b) agrees with an important
 calculation of Deligne's in [D], 10.1, 
which we now describe: Using the $C^{\infty}$ splitting of $\E$, we
can write the $\ov{\partial}$ operator for $E$ in matrix form:
\[
\dis
\ov{\partial}_E=\left(
\begin{array}{cc}
\ov{\partial}_S & \alpha \\
0 & \ov{\partial}_Q
\end{array}
\right),
\ \ \ \mbox{for some } \alpha\hin A^{0,1}(X,\Hom(Q,S)).
\]
Let $\alpha^*\hin A^{1,0}(X,\Hom(S,Q))$ be the transpose of $\alpha$,
defined using complex conjugation of forms and the metrics $h_S$ and
$h_Q$. If $\nabla$ is the induced connection on $\Hom(Q,S)$, we can write
\[
\dis
K_E=\left(
\begin{array}{c|c}
K_S-\frac{i}{2\pi}\alpha\alpha^* & \nabla^{1,0}\alpha \\ \hline
-\nabla^{0,1}\alpha^* & K_Q-\frac{i}{2\pi}\alpha^*\alpha 
\end{array}
\right).
\]
 Thus Corollary \ref{c1}(b) implies that
\[
\wt{c_2}(\ov{\E})=-\frac{1}{2\pi i}\Tr (\alpha\alpha^*)=
\frac{1}{2\pi i}\Tr (\alpha^*\alpha),
\]
and we have recovered Deligne's result. In this form the
 calculation was used by A. Moriwaki and
C. Soul\'{e} to obtain a Bogomolov-Gieseker type inequality and a
Kodaira vanishing theorem on 
arithmetic surfaces, respectively (see [Mo] and [S]).

The calculation of $\wt{c_2}$ shows that in general Bott-Chern
forms are not closed. In fact, calculating $\wt{c_k}$ for $k \gequ 3$
leads to much more complicated formulas, involving traces of products
of curvature matrices, for which a clear geometric
interpretation is lacking 
(unlike the matrix $\alpha$ above, whose negative transpose
 $-\alpha^*$ is the second fundamental form of $\ov{\E}$). 
In the next two sections we shall see that when 
 $\ov{E}$ is a projectively flat bundle, the Bott-Chern forms are 
closed and can be calculated explicitly for any $\phi\hin I(n)$.

\section{$0\ra \ov{S} \ra \ov{E}\ra \ov{Q}\ra 0$ with $\ov{E}$ flat}
\label{flatses}

Throughout this section we will assume that the hermitian vector 
bundle $\ov{E}$ is {\em flat}, i.e. that $K_E=0$. As before,
the metrics $h_S$ and $h_Q$ will be induced from the metric
on $E$. Define the {\em harmonic numbers}
$\dis \H_k=\sum_{i=1}^k\frac{1}{i}$, $\H_0=0$. 

Let $\lambda$ be a partition of
$k$ (we denote this by $\lambda\vdash k$). Recall that the polynomials
$\{p_{\lambda} : \lambda\vdash k\}$ form a $\Q$-basis for the vector
space of symmetric homogeneous polynomials in $x_1,\ldots,x_n$ of 
degree $k$. The following result
 computes the Bott-Chern form corresponding to any such invariant
polynomial:

\begin{thm} \label{bcflat}
The Bott-Chern class $\wt{p_{\lambda}}(\ov{\E})$ in $\wt{A}(X)$ is
the class of

\ \ \ \ \ \ \ \ \ \
{\em (i)} $k\H_{k-1}p_{k-1}(\ov{Q})$, if $\lambda=k=(k,0,0,\ldots,0)$

\ \ \ \ \ \ \ \ \ \
{\em (ii)} 0, otherwise.

\end{thm}
{\bf Proof.} Let us first compute $\wt{p_k}(\ov{\E})$ for $p_k(A)=
\Tr (A^k)$. Since $K_E=0$, the deformed matrix 
$K(u)=(1-u)K_{S\oplus Q}$, where $K_{S\oplus Q}=
\left(
\begin{array}{c|c}
K_S & 0 \\ \hline
0 & K_Q
\end{array} \right)$. Since
\[
\dis
\int_0^1\frac{(1-u)^{k-1}-1}{u}\,du=-\int_0^1\frac{t^{k-1}-1}{t-1}\,dt=
-\H_{k-1},
\]
we obtain
\[
\dis
\wt{p_k}(\ov{\E})=-\H_{k-1}p_k^{\pr}(K_{S\oplus Q}; J_r)=
-k\H_{k-1}\Tr (K_S^{k-1})=-k\H_{k-1}p_{k-1}(\ov{S}).
\]
Now since $p_k(\ov{S}\bigoplus\ov{Q})-p_k(\ov{E})$ is exact, 
$p_k(\ov{E})=0$ and $p_k(\ov{S}\bigoplus\ov{Q})=
p_k(\ov{S})+p_k(\ov{Q})$, we conclude that $p_k(\ov{S})=-p_k(\ov{Q})$
in $\wt{A}(X)$, for each $k \gequ 1$. This proves (i).

Let $\lambda=(\lambda_1,\lambda_2,\ldots,\lambda_m)$ be a partition 
($m\gequ 2$). 
Proposition \ref{bcprop} implies that
\[
\dis
\wt{p_{\lambda}}(\ov{\E})=\wt{p_{\lambda_1}}(\ov{\E})
p_{\lambda_2}\cdots p_{\lambda_m}(\ov{S}\bigoplus\ov{Q}).
\]
But $\wt{p_{\lambda_1}}(\ov{\E})$ is a  closed form (by (i)), and
$p_{\lambda_2}\cdots p_{\lambda_m}(\ov{S}\bigoplus\ov{Q})$ is an
exact form. Thus $\wt{p_{\lambda}}(\ov{\E})$ is exact, and so vanishes
in $\wt{A}(X)$.\
\endproof

\medskip

It follows from Theorem \ref{bcflat} that for any $\phi \hin I(n)$, the
Bott-Chern form $\wt{\phi}(\ov{\E})$ is a linear combination of homogeneous
components of the Chern character form $ch(\ov{Q})$.
In [Ma], Theorem 3.4.1 we find the calculation
\begin{equation}
\label{maillotcalc}
\wt{c_k}(\ov{\E})=\H_{k-1}\sum_{i=0}^{k-1}ic_i(\ov{S})c_{k-1-i}(\ov{Q})
\end{equation}
for the Chern forms $\wt{c_k}$. Our result gives the following

\begin{prop} \label{maillotprop}
$\wt{c_k}(\ov{\E})=(-1)^{k-1}\H_{k-1}p_{k-1}(\ov{Q})$.
\end{prop}
{\bf Proof.} By Newton's identity (\ref{ni}) we have
\begin{equation}
\label{newton}
\wt{p_k}-\wt{c_1p_{k-1}}+\wt{c_2p_{k-2}}-\cdots+(-1)^kk\wt{c_k}=0.
\end{equation}

Reasoning as in Theorem \ref{bcflat}, we see that if $\phi$ and $\psi$ are
two homogeneous invariant polynomials of positive degree, then
$\wt{\phi\psi}(\ov{\E})=0$ in $\wt{A}(X)$. Thus (\ref{newton})
gives 
$\dis
\wt{c_k}(\ov{\E})=\frac{(-1)^{k-1}}{k}\wt{p_k}(\ov{\E})=
(-1)^{k-1}\H_{k-1}p_{k-1}(\ov{Q})$. \endproof

\medskip

\noindent
{\bf Remark.} The result of Proposition \ref{maillotprop} agrees with
(\ref{maillotcalc}), i.e.
 $\dis (-1)^kp_k(\ov{Q})=\sum_{i=0}^k
ic_i(\ov{S})c_{k-i}(\ov{Q})$ in $\wt{A}(X)$. To see this, let
$h(t)=\sum c_i(\ov{S})t^i$, $g(t)=\sum c_j(\ov{Q})t^j$, and
$f(t)=\sum ic_i(\ov{S})t^i$. Then $h(t)g(t)=1$ in $\wt{A}(X)[t]$, and
$f(t)=th^{\pr}(t)$. Choose formal variables 
$\{x_{\alpha}\}_{1\lequ\alpha\lequ r}$ and set $c_i(\ov{S})=
e_i(x_1,\ldots,x_r)$, so that $\dis h(t)=\prod_{\alpha}(1+x_{\alpha}t)$.
Then $\dis f(t)=\sum_{\alpha}tx_{\alpha}\prod_{\beta\neq\alpha}
(1+x_{\beta}t)$. Thus
\[
\dis
f(t)g(t)=\frac{f(t)}{h(t)}=\sum_{\alpha}\frac{x_{\alpha}t}{1+x_{\alpha}t}=
r-\sum_{\alpha}\frac{1}{1+x_{\alpha}t}=
\]
\[
\dis
r-\sum_{\alpha,i}(-1)^ix_{\alpha}^it^i=
r-\sum_i(-1)^ip_i(\ov{S})t^i=r+\sum_i(-1)^ip_i(\ov{Q})t^i.
\]
Comparing coefficients of $t^k$ on both sides gives the result.

\medskip

 We can use Theorem \ref{bcflat} to calculate $\wt{\phi}
(\ov{\E})$ for $\phi\hin I(n)_k$: it is enough to find the
coefficient of the power sum $p_k$ when $\phi$ is expressed as a
linear combination of $\{p_{\lambda}\}_{\lambda\vdash k}$ in $\L(n,\Q)$.
 For example, we have

\begin{cor} \label{appl0}
$\dis \wt{ch}(\ov{\E})=\sum_k\H_k ch_k(\ov{Q})$, where $ch_k$ denotes
the $k$-th homogeneous component of the Chern character form.
\end{cor}

\begin{cor} \label{appl1}
Let $\lambda$ be a partition of $k$ and $s_{\lambda}$ the
corresponding Schur polynomial in $\L(n,\Q)$.
Then $\wt{s_{\lambda}}(\ov{\E})=0$ unless $\lambda$ is a hook
$\lambda_i=(i,1,1,\ldots,1)$, in which case
$\wt{s_{\lambda_i}}(\ov{\E})=(-1)^{k-i}\H_{k-1}p_{k-1}(\ov{Q})$.
\end{cor}
{\bf Proof.} The proof is based on the Frobenius formula
\[
\dis
s_{\lambda}=\frac{1}{k!}\sum_{\sigma\in S_k}\chi_{\lambda}(\sigma)
p_{(\sigma)}
\]
where $(\sigma)$ denotes the partition of $k$ determined by the cycle 
structure of $\sigma$ (cf. [M], \S I.7).
By the above remark, 
$\wt{s_{\lambda}}(\ov{\E})=\chi_{\lambda}((12\ldots k))\H_{k-1}
p_{k-1}(\ov{Q})$. Using the combinatorial rule for computing 
$\chi_{\lambda}$ found in [M], p. 117, Example 5, we obtain
\[
\dis
\chi_{\lambda}((12\ldots k))=
\left\{ \begin{array}{cl}
(-1)^{k-i}, & \mbox{if } \lambda=\lambda_i \mbox{ is a hook} \\
0, & \mbox{otherwise}. 
\end{array}
\right.
\]
\endproof

The most natural instance of a sequence $\ov{\E}$ with $\ov{E}$
flat is the classifying sequence over the grassmannian $G(r,n)$.
As we shall see in \S \ref{grass}, the calculation of 
Bott-Chern forms for this sequence leads to a presentation of
arithmetic intersection ring of the arithmetic grassmannian over $\Spec\Z$.

\section{Calculations when $\ov{E}$ is projectively flat}
\label{projflat}

We will now generalize the results of the last section to the case 
where $E$ is {\em projectively flat}, i.e. the curvature matrix
$K_E$ of $\ov{E}$ is a multiple of the identity matrix:
$K_E=\omega Id_n$. This is true if
$\dis E=\ov{L}^{\oplus n}$ for some hermitian line bundle $\ov{L}$,
with $\omega=c_1(\ov{L})$ the first Chern form of $\ov{L}$.

The Bott-Chern forms (for the induced metrics) 
are always closed in this case as well, and 
will be expressed in terms of characteristic classes of the bundles
involved. However this seems to be the most general case where this phenomenon
occurs.

The key observation is that for projectively flat bundles, the curvature
matrix $K_E=\omega Id_n$ in {\em any} local trivialization. Thus we have
\[
\dis
K(u)=
\left(
\begin{array}{c|c}
(1-u)K_S+u\omega Id_r & 0 \\ \hline
0 & (1-u)K_Q+u\omega Id_s
\end{array} \right)
\]
where $s=n-r$ denotes the rank of $Q$. Now Theorem
\ref{calc} gives
\[
\dis
\wt{p_k}(\ov{\E})=
k\int_0^1\frac{1}{u}\Tr[(u\omega Id_r+(1-u)K_S)^{k-1}-K_S^{k-1}]\,du=
\]
\[
-k\H_{k-1}p_{k-1}(\ov{S})+k\sum_{j=1}^{k-1}{k-1 \choose j}
\Tr(\omega^jK_S^{k-1-j})
\int_0^1u^{j-1}(1-u)^{k-j-1}\,du.
\]
Integrating by parts gives $\dis
\int_0^1u^m(1-u)^n\,du=\frac{1}{m+n+1}{m+n \choose n}^{-1}$, thus
\begin{equation}
\label{pkS}
\frac{1}{k}\wt{p_k}(\ov{\E})=
-\H_{k-1}p_{k-1}(\ov{S})+\sum_{j=1}^{k-1}\frac{\omega^j}{j}
p_{k-1-j}(\ov{S}).
\end{equation}

 We can rewrite this as an equation involving power sums of the
quotient bundle: since $p_k(\ov{S})+p_k(\ov{Q})-p_k(\ov{L}^{\oplus n})=0$
in $\wt{A}(X)$, we have $p_k(\ov{S})=n\omega^k-p_k(\ov{Q})$. Thus
(\ref{pkS}) becomes
\begin{equation}
\label{psum}
\frac{1}{k}\wt{p_k}(\ov{\E})=
\H_{k-1}p_{k-1}(\ov{Q})-\sum_{j=1}^{k-1}\frac{\omega^j}{j}
p_{k-1-j}(\ov{Q}).
\end{equation}
 
\begin{thm} \label{appl2}
Let $X$ be a complex manifold, $\ov{E}$ a projectively flat
hermitian vector bundle
over $X$. Let $0\ra \ov{S}\ra\ov{E}\ra\ov{Q}\ra 0$ a short exact 
sequence of vector bundles over $X$ with metrics on $S$, $Q$ induced from 
$\ov{E}$. Then for any invariant polynomial $\phi\hin I(n)$,
$\phi(\ov{S}\bigoplus\ov{Q})=\phi(\ov{E})$ as differential forms
on $X$.
\end{thm}
{\bf Proof.} Since the
$p_{\lambda}$ form an additive basis for $I(n)$, it suffices to prove the
result when $\phi=p_{\lambda}$. 
The above calculation shows that $\wt{p_k}$ is a closed
form. This combined with Proposition \ref{bcprop} shows that
 $\wt{p_{\lambda}}$ is closed for any partition $\lambda$. 
Thus
\[
\dis
p_{\lambda}(\ov{S}\bigoplus\ov{Q})-p_{\lambda}(\ov{E})=
dd^c\wt{p_{\lambda}}=0. 
\] \endproof

\noindent
{\bf Remark.}
If $E$ is a trivial vector bundle, this result follows by pulling
back the exact sequence $\ov{\E}$ from the classifying sequence on the
Grassmannian. The forms are equal there because they are invariant 
with respect to the $U(n)$ action, so harmonic.

\medskip

  The Bott-Chern forms $\wt{p_{\lambda}}$ for a general partition
$\lambda=(\lambda_1,\ldots,\lambda_m)$ can be computed by using 
Proposition \ref{bcprop}. If $|\lambda|=\sum\lambda_i=k$ then we have
\begin{equation}
\label{ppsum}
\dis
\wt{p_{\lambda}}(\ov{\E})=\sum_{i=1}^m
\wt{p_{\lambda_i}}(\ov{\E}) \prod_{j\neq i}
p_{\lambda_j}(\ov{E})=
n^{m-1}\sum_{i=1}^m\omega^{k-\lambda_i}\wt{p_{\lambda_i}}(\ov{\E}).
\end{equation}
In principle equations (\ref{psum}) and (\ref{ppsum})
 can be used to compute $\wt{\phi}(\ov{\E})$
for any characteristic class $\phi$.

  We now find a more explicit formula for
 the Bott-Chern forms of Chern classes. The 
computation is not as straightforward, as the argument of 
Proposition \ref{maillotprop} does not apply. 
Since by Theorem \ref{calc} the calculation depends only on the
 curvature matrices $K_E$, $K_S$ and $K_Q$, we may assume
\[
\dis
\ov{\E}\ : \ \
   0\ra \ov{S} \ra \ov{L}\otimes\C^n\ra \ov{Q}\ra 0
\]
is our chosen sequence, and define a new sequence
\[
\dis
\ov{\E^{\pr}}=\ov{\E}\otimes\ov{L}^*\ : \ 
   0\ra \ov{S}\otimes\ov{L}^* \ra\C^n\ra \ov{Q}\otimes\ov{L}^* \ra 0.
\]

The metrics on the bundles in $\ov{\E^{\pr}}$ are induced from the
trivial metric on $\C^n$. Using Propositions \ref{twist} and 
\ref{maillotprop} now gives
\[
\dis
\wt{c_k}(\ov{\E})=
\wt{c_k}(\ov{\E^{\pr}}\otimes\ov{L})=
\sum_{i=1}^k{n-i \choose k-i}\wt{c_i}(\ov{\E^{\pr}})c_1(\ov{L})^{k-i}=
\]
\[
\dis
\sum_{i=1}^k{n-i \choose k-i}(-1)^{i-1}\H_{i-1}
p_{i-1}(\ov{Q}\otimes\ov{L}^*)\omega^{k-i}=
\]
\[
\dis
\sum_{i=1}^k\sum_{j=0}^{i-1}(-1)^j{n-i \choose k-i}{i-1 \choose j}\H_{i-1}
\omega^{k-1-j}p_j(\ov{Q})=
\]
\[
\dis
\sum_{j=0}^{k-1}(-1)^jd_j\omega^{k-1-j}p_j(\ov{Q}),
\]
where
\[
\dis
d_j=\sum_{i=j+1}^k{n-i \choose k-i}{i-1 \choose j}\H_{i-1}.
\]

 To find a closed form for the sum $d_j$, we can use the general identity
\begin{equation}
\label{harm}
\sum_{i=q-s}^{n-p}{n-i \choose p}{s+i \choose q}\H_{s+i}=
{n+s+1 \choose p+q+1}(\H_{n+s+1}-\H_{p+q+1}+\H_p).
\end{equation}
This is identity (10) in [Sp]. In passing we note that writing
equation (\ref{harm}) without the harmonic number terms:
\[
\dis
\sum_{i=q-s}^{n-p}{n-i \choose p}{s+i \choose q}={n+s+1 \choose p+q+1}
\]
gives a well known identity among binomial coefficients.
Applying (\ref{harm}) to $d_j$ and replacing
$k$ by $k+1$ and $j$ by $k-i$ we arrive at the formula
\[
\dis
\wt{c_{k+1}}(\ov{\E})=
\sum_{i=0}^k (-1)^{k-i} {n \choose i} {\frak H}_i \omega^i p_{k-i}(\ov{Q}),
\]
where ${\frak H}_i=\H_n-\H_{n-i}+\H_{k-i}$. As remarked previously,
this calculation is 
valid for any projectively flat bundle $\ov{E}$ with $c_1(\ov{E})=n\omega$.

 Of course one can use the above method to compute the Bott-Chern form
$\wt{p_k}(\ov{\E})$ as well; however this leads to a more complicated formula
than (\ref{psum}). Equating the two proves(!) the following interesting
combinatorial identity (compare [Sp], identity (30)):

\begin{equation}
\label{identity}
\sum_{i=0}^s(-1)^{i+1} {n \choose i,s-i} \H_{n-s+i}=\frac{1}{s}\ \ \ \ \ \ \ \
(n\gequ s).
\end{equation}
Here $\dis {n \choose i,j}$ is a trinomial coefficient. 

The following summarizes the calculations of this section:

\begin{thm}
\label{bccalcthm}
Let $X$ be a complex manifold, $\ov{E}$ a projectively flat 
hermitian vector bundle over $X$, with $c_1(\ov{E})=n\omega$.
 Let $0\ra \ov{S}\ra\ov{E}\ra\ov{Q}\ra 0$ be a short exact 
sequence of vector bundles over $X$ with metrics on $S$, $Q$ induced from 
$\ov{E}$. Then
\[
\dis
\wt{p_{k+1}}(\ov{\E})=
(k+1)\H_kp_k(\ov{Q})-(k+1)\sum_{i=1}^k \frac{\omega^i}{i} p_{k-i}(\ov{Q})
\]
\[
\dis
\wt{c_{k+1}}(\ov{\E})=
\sum_{i=0}^k (-1)^{k-i} {n \choose i} {\frak H}_i \omega^i p_{k-i}(\ov{Q})
\]
where ${\frak H}_i=\H_n-\H_{n-i}+\H_{k-i}$.
\end{thm}

Note that the formulas in Theorem \ref{bccalcthm}
reduce to the ones of the previous section
when $\omega=0$!

\section{Arithmetic intersection theory}
\label{ait}

 We recall here the generalization of Arakelov
theory to higher dimensions due to H. Gillet and C. Soul\'{e}. 
Our main references are [GS1], [GS2] and the exposition in [SABK].
 For $A$ an abelian group, $A_{\Q}$ denotes $A\otimes_{\Z}\Q$.
Let $X$ be an
{\em arithmetic scheme over $\Z$}, by which we mean a regular scheme, 
projective and flat over $\mbox{Spec}\Z$.
For $p\gequ 0$, let $X^{(p)}$ be the set of integral subschemes of $X$ of
codimension $p$ and 
 $Z^p(X)$ be the group 
of codimension $p$ cycles on $X$. The $p$-th Chow group of $X$:
$CH^p(X):=Z^p(X)/R^p(X)$, where $R^p(X)$ is the subgroup of $Z^p(X)$ 
generated by the cycles $\mbox{div}f$, $f\hin k(x)^*$, $x\hin X^{(p-1)}$.
Let $CH(X)=\bigoplus_pCH^p(X)$. If $X$ is smooth over $\Spec\Z$, then
the methods of [F] can be used to give $CH(X)$ the structure of a 
commutative ring. In general one has a product structure on $CH(X)_{\Q}$
after tensoring with $\Q$.

Let $D^{p,p}(X(\C))$ denote the space of complex currents of type $(p,p)$
on $X(\C)$, and $F_{\infty}:X(\C)\ra X(\C)$ the involution induced by
complex conjugation. Let $D^{p,p}(X_{\R})$ (resp. $A^{p,p}(X_{\R})$)
be the subspace of $D^{p,p}(X(\C))$ (resp. $A^{p,p}(X(\C))$) generated
by real currents (resp. forms) $T$ such that $F^*_{\infty}T=(-1)^pT$;
denote by $\wt{D}^{p,p}(X_{\R})$ and $\wt{A}^{p,p}(X_{\R})$
the respective images in $\wt{D}^{p,p}(X(\C))$ and
$\wt{A}^{p,p}(X(\C))$. 

An {\em arithmetic cycle} on $X$ of codimension $p$ is a pair $(Z,g_Z)$ in 
the group
$Z^p(X)\bigoplus\wt{D}^{p-1,p-1}(X_{\R})$, where $g_Z$ is a 
{\em Green current} for $Z(\C)$, i.e. a current such that 
$dd^cg_Z+\delta_{Z(\C)}$ is represented by a smooth form.
 The group of arithmetic cycles is
denoted by $\wh{Z}^p(X)$. If $x\hin X^{(p-1)}$ and $f\hin
k(x)^*$, we let $\wh{\mbox{div}}f$ denote the arithmetic cycle
$(\mbox{div}f, [-\log|f_{\C}|^2\cdot \delta_{x(\C)}])$.

The {\em $p$-th arithmetic Chow group of $X$}: $\wh{CH}^p(X):=
\wh{Z}^p(X)/\wh{R}^p(X)$, where $\wh{R}^p(X)$ is the
subgroup of $\wh{Z}^p(X)$ generated by the cycles
$\wh{\mbox{div}}f$, $f\hin k(x)^*$, $x\hin X^{(p-1)}$.
Let $\wh{CH}(X)=\bigoplus_p\wh{CH}^p(X)$.

 We have the following canonical morphisms of abelian groups:
\[
\dis
\zeta :\wh{CH}^p(X)  \longrightarrow  CH^p(X), \ \ \ 
 {[(Z,g_Z)]} \longmapsto  {[Z]},
\]
\[
\dis
 \omega : \wh{CH}^p(X) \longrightarrow  \Ker d\cap
 \Ker d^c\hspace{0.2cm} (\subset A^{p,p}(X_{\R})), \ \ \
 {[(Z,g_Z)]}  \longmapsto dd^cg_Z+\delta_{Z(\C)},
\]
\[
\dis
 a : \wt{A}^{p-1,p-1}(X_{\R}) \longrightarrow \wh{CH}^p(X), \ \ \
 \eta  \longmapsto  {[(0,\eta)]}.
\]

One can define a pairing $\wh{CH}^p(X)\otimes\wh{CH}^q(X)
\ra \wh{CH}^{p+q}(X)_{\Q}$ which turns $\wh{CH}(X)_{\Q}$
into a commutative graded unitary $\Q$-algebra. The maps $\zeta$, $\omega$
are $\Q$-algebra homomorphisms. If $X$ is smooth over $\Z$
one does not have to tensor with $\Q$.
The definition of this pairing is difficult; the construction uses the
{\em star product} of Green currents, which in turn relies upon Hironaka's
resolution of singularities to get to the case of divisors.
The functor $\wh{CH}^p(X)$ is 
contravariant in $X$, and covariant for proper maps which are smooth on
the generic fiber.
 
Choose a 
 K\"{a}hler form $\omega_0$ on $X(\C)$ such that $F^*_{\infty}\omega_0=
-\omega_0$ (this is equivalent to requiring that the corresponding K\"ahler
metric is invariant under $F_{\infty}$).
 It is natural to utilize the theory of harmonic forms on 
$X$ in the study of Green currents on $X(\C)$. Following [GS1], we call
the pair $\ov{X}=(X,\omega_0)$ an {\em Arakelov variety}. 
By the Hodge decomposition
theorem, we have $A^{p,p}(X_{\R})=\H^{p,p}(X_{\R})\oplus \mbox{Im}d\oplus
\mbox{Im}d^*$, where $\H^{p,p}(X_{\R})=\Ker\Delta\subset A^{p,p}(X)$
denotes the space of harmonic (with respect to $\omega_0$) $(p,p)$ 
forms $\alpha$ on $X(\C)$ such that $F_{\infty}^*\alpha=(-1)^p\alpha$.
The subgroup $CH^p(\ov{X}):=\omega^{-1}(\H^{p,p}(X_{\R}))$ of 
$\wh{CH}^p(X)$ is called the {\em $p$-th Arakelov Chow group of $X$}.
Let $CH(\ov{X})=\bigoplus_{p\gequ 0} CH^p(\ov{X})$.
 $CH^p(\ov{X})$ is a direct summand of $\wh{CH}^p(X)$, and there is
an exact sequence 
\begin{equation}
\label{ex1}
CH^{p,p-1}(X) \stackrel{\rho}\longrightarrow \H^{p-1,p-1}(X_{\R})
\stackrel{a}\longrightarrow CH^p(\ov{X}) 
\stackrel{\zeta}\longrightarrow CH^p(X)\longrightarrow 0.
\end{equation}
In the above sequence the group $CH^{p,p-1}(X)$ is defined as the 
$E_2^{p,1-p}$ term of a certain spectral sequence used by Quillen to 
calculate the higher algebraic $K$-theory of $X$, and the map $\rho$ 
coincides with the Beilinson regulator map (cf. [G] and [GS1], 3.5).

If $\H(X_{\R})=\bigoplus_p \H^{p,p}(X_{\R})$ is a subring of
$\bigoplus_p A^{p,p}(X_{\R})$, for example if $X(\C)$ is a 
curve, an abelian variety or a hermitian symmetric space (e.g. a
grassmannian), then $CH(\ov{X})_{\Q}$ is a subring of 
$\wh{CH}(X)_{\Q}$.  This is not the case in general; for
example it fails to be true for the complete flag varieties.

Arakelov [A] introduced the group $CH^1(\ov{X})$, where $\ov{X}=(X,g_0)$
is an arithmetic surface with the metric $g_0$ on the Riemann surface
$X(\C)$ given by $ \frac{i}{2g}\sum \omega_j\wedge \ov{\omega}_j$.
Here $g$ is the genus of $X(\C)$ and $\{\omega_j\}$ for $1\lequ j \lequ g$
is an orthonormal basis of the space of holomorphic one forms on $X(\C)$.

  A {\em hermitian vector bundle} $\ov{E}=(E,h)$ on an arithmetic scheme $X$ is
an algebraic vector bundle $E$ on $X$ such that the induced holomorphic
vector bundle $E(\C)$ on $X(\C)$ has a hermitian metric $h$, which is 
invariant under complex conjugation, i.e. $F_{\infty}^*(h)=h$.

 To any hermitian vector bundle one can attach characteristic classes 
$\wh{\phi}(\ov{E})\hin \wh{CH}(X)_{\Q}$, for any $\phi\hin I(n,\Q)$, 
where $n=\rk E$. For example, we have {\em arithmetic Chern classes} 
$\wh{c}_k(\ov{E}) \hin \wh{CH}^k(X)$. 
Some basic properties of these classes are:

\noindent 
(1) $\wh{c}_0(\ov{E})=1$ and $\wh{c}_k(\ov{E})=0$ for $k>\mbox{rk}E$.

\noindent
(2) The form $\omega(\wh{c}_k(\ov{E}))=c_k(\ov{E})\hin A^{k,k}(X_{\R})$
is the $k$-th Chern form of the hermitian bundle $\ov{E(\C)}$.

\noindent
(3) $\zeta(\wh{c}_k(\ov{E}))=c_k(E)\hin CH^k(X)$.

\noindent
(4) $f^*\wh{c}_k(\ov{E})=\wh{c}_k(f^*\ov{E})$, for every morphism
$f:X\ra Y$ of regular schemes, projective and flat over $\Z$.

\noindent
(5) If $\ov{L}$ is a hermitian line bundle, $\wh{c}_1(\ov{L})$ is the
class of $(\mbox{div}(s),-\log \|s\|^2)$ for any rational section $s$
of $L$.

 Analogous properties are satisfied by $\wh{\phi}$ for any $\phi\hin
I(n,\Q)$ (see [GS2], Th. 4.1). 
The most relevant property of these
characteristic classes is their behaviour in short exact sequences: 
if 
\[
\dis
\ov{\E} :\ 0 \ra \ov{S} \ra \ov{E} \ra \ov{Q} \ra 0
\]
is such a sequence of hermitian vector bundles over $X$, then
\begin{equation}
\label{key}
\wh{\phi}(\ov{S}\oplus\ov{Q})-
\wh{\phi}(\ov{E})=a(\wt{\phi}(\ov{\E})). 
\end{equation}
Relation (\ref{key}) is the main tool for calculating
intersection products of classes in $\wh{CH}(X)$ that come from 
characteristic classes of vector bundles.
Combining it with the results of \S 4 and \S 5 gives immediate consequences for
such intersections. For example, we have

\begin{cor} \label{appl3}
Let
$\ov{\E}: 0\ra \ov{S}\ra\ov{E}\ra\ov{Q}\ra 0$ be a short exact 
sequence of hermitian vector bundles over an arithmetic scheme $X$.
 Assume that the metrics on $S(\C)$, $Q(\C)$ are 
induced from that on 
$E(\C)$.

\noindent
(a) If $\ov{E(\C)}$ is flat, then
\[
\dis
(1)\ \ \ \wh{p_{\lambda}}(\ov{S}\oplus\ov{Q})=\wh{p_{\lambda}}(\ov{E}),
\ \mbox{ if } \lambda \mbox{ has length } >1, \mbox{ and}
\]
\[
\dis
(2)\ \ \ \wh{p_k}(\ov{S})+\wh{p_k}(\ov{Q})-\wh{p_k}(\ov{E})=
k\H_{k-1}a(p_{k-1}(\ov{Q})), \ \forall k\gequ 1,
\]
in the arithmetic Chow group $\wh{CH}(X)_{\Q}$.

\noindent
(b) If $\ov{E}=\ov{L}^{\oplus n}$ for some hermitian line bundle $\ov{L}$
and $\omega=c_1(\ov{L(\C)})$, then
\[
\dis
\wh{c}(\ov{S})\wh{c}(\ov{Q})-\wh{c}(\ov{E})=
\sum_{i,j} (-1)^j {n \choose i} (\H_n-\H_{n-i}+\H_j) a(\omega^i p_j(\ov{Q})),
\]
in the arithmetic Chow group $\wh{CH}(X)$.
 
\end{cor}

\section{Arakelov Chow rings of grassmannians} 
\label{grass}

In this section $G=G(r,n)$ will denote the grassmannian over $\Spec\Z$.
Over any field $k$, $G$
parametrizes the $r$-dimensional linear subspaces
of a vector space over $k$. 
Let 
\begin{equation}
\label{mex}
\ov{\E}:\ 0\ra \ov{S} \ra \ov{E}\ra \ov{Q} \ra 0
\end{equation}
denote the universal exact sequence of vector bundles over $G$. Here 
the trivial bundle $E(\C)$ is given the trivial metric and the tautological
subbundle $S(\C)$ and quotient bundle $Q(\C)$ the induced metrics. 
The homogeneous space $G(\C)\simeq U(n)/(U(r)\times U(n-r))$ 
is a complex manifold. $G(\C)$ is endowed 
with a natural $U(n)$-invariant metric coming from the K\"{a}hler form
$\eta_G=c_1(\ov{Q(\C)})$.

$G$ is a smooth arithmetic scheme and $G(\C)$ with the metric coming
from $\eta_G$ is a hermitian symmetric
space, so we have an Arakelov Chow ring $CH(\ov{G})$.
Note that since the hermitian vector bundles in (\ref{mex})
are invariant under the action of $U(n)$, their Chern
forms are harmonic, and thus the arithmetic characteristic classes
obtained are all elements of $CH(\ov{G})$.
V. Maillot [Ma] found a presentation of 
$CH(\ov{G})$, using the above observation and the short exact sequence
(\ref{ex1}). We wish to offer another description
of this ring, based on the calculations in this paper.

First recall the geometric picture: for the ordinary Chow ring we have
\begin{equation}
\label{chowg}
CH(G)=\frac{\Z[c(S),c(Q)]}{\left<c(S)c(Q)=1\right>}.
\end{equation}
(see for instance [F], Example 14.6.6).
If $x_1,\ldots,x_r$ are the Chern roots of $S$,
$y_1,\ldots,y_s$ are the Chern roots of $Q$, $H=S_r\times S_{n-r}$ is
the product of two symmetric groups, and
$t$ is a formal variable, 
then (\ref{chowg}) can be
rewritten 
\begin{equation}
\label{chow}
CH(G)=\frac{\Z[x_1,\ldots,x_r,y_1,\ldots,y_s]^H}
{\left<\prod_i(1+x_it)\prod_j(1+y_jt)=1\right>}.
\end{equation}

 Maillot's presentation of $CH(\ov{G})$
is an analogue of (\ref{chowg}); ours
will be an analogue of (\ref{chow}). We introduce $2n$ variables
\[
\wh{x}_1,\ldots,\wh{x}_r,\wh{y}_1,\ldots,\wh{y}_s,
x_1,\ldots,x_r,y_1,\ldots,y_s
\]
 and consider the rings
\[
A=\Z[\wh{x}_1,\ldots,\wh{x}_r,\wh{y}_1,\ldots,\wh{y}_s]^H\ \ 
\mbox{ and } \ \
B=\R[x_1,\ldots,x_r,y_1,\ldots,y_s]^H
\]
 and the ring homomorphism
$\omega:A\ra B$ defined by $\omega(\wh{x}_i)=x_i$ and
$\omega(\wh{y}_j)=y_j$. A ring structure is defined on the 
abelian group $A\oplus B$ by setting
\[
\dis
(\wh{x},x^{\pr})*(\wh{y},y^{\pr})=
(\wh{x}\wh{y},\omega(\wh{x})y^{\pr}+
x^{\pr}\omega(\wh{y})).
\]
We will adopt the convention that $\wh{\alpha}$ denotes $(\wh{\alpha},0)$,
$\beta$ denotes $(0,\beta)$, 
and any product $\prod x_iy_j$ denotes
$(0,\prod x_iy_j)$; the multiplication $*$ is thus characterized by
the properties $\wh{\alpha} * \beta=\alpha\beta$ and $\beta_1 * \beta_2=0$.
We now define two sets of relations in $(A\oplus B)[t]$:
\[
\dis
{\cal R}_1:\
\prod_i(1+x_it)\prod_j(1+y_jt)=1,
\]
\[
\dis
{\cal R}_2:\ 
\prod_i(1+\wh{x}_it) * \prod_j(1+\wh{y}_jt) *
\left(1+t\sum_j\frac{\log(1+y_jt)}{1+y_jt}\right)=1.
\]
and let ${\cal A}$ denote the quotient of the graded ring
$A\oplus B$ by these relations. Using this notation we can state

\begin{thm}
\label{presentation}
There is a unique ring isomorphism
$\Phi: {\cal A}\ra CH(\ov{G})$ such that
\[
\dis
\Phi(\prod_i (1+\wh{x}_it^i))=
\sum_i \wh{c_i}(\ov{S})t^i,\ \ \ \ 
\Phi(\prod_j (1+\wh{y}_jt^j))=
\sum_j \wh{c_j}(\ov{Q})t^j,
\]
\[
\dis
\Phi(\prod_i (1+x_it^i))=
\sum_i a(c_i(\ov{S}))t^i,\ \ \ \
\Phi(\prod_j (1+y_jt^j))=
\sum_j a(c_j(\ov{Q}))t^j.
\]
\end{thm}

\medskip
{\bf Proof.} The isomorphism $\Phi$ of ${\cal A}$ with $CH(\ov{G})$ is 
obtained exactly as in [Ma], Theorem 4.0.5. The key fact is that since $G$
has a cellular decomposition (in the sense of [F], Ex. 1.9.1.), it follows
that $CH^{p,p-1}(G)=0$ for all $p$ (using the excision exact sequence
for the groups $CH^{*,*}(G)$; cf. [G], \S 8). Summing the sequence
(\ref{ex1}) over all $p$ gives
\begin{equation}
\label{equ}
0 \longrightarrow \H(G_{\R})
\stackrel{a}\longrightarrow CH(\ov{G}) 
\stackrel{\zeta}\longrightarrow CH(G)\longrightarrow 0.
\end{equation}
We can now use our knowledge of the rings $\H(G_{\R})$
and $CH(G)$ together with the five lemma, as in loc. cit.
The multiplication $*$ is a consequence of the general
identity $a(x)y=a(x\omega(y))$ in $\wh{CH}(G)$.
To complete the argument we
must show that the relation $\wh{c}(\ov{S})\wh{c}(\ov{Q})=1+
a(\wt{c}(\ov{\E}))$ translates to the relation ${\cal R}_2$ above.

Let $p_i(y)$ be the $i$-th power sum in the variables
$y_1\ldots,y_s$, identified under $\Phi$ with the 
class $a(p_i(\ov{Q}))$ in $CH(\ov{G})$ (we will use such 
identifications freely in the sequel).
We also define $\dis p_a(t)=\sum_{i=0}^{\infty}
(-1)^{i+1}\H_ip_i(y)t^{i+1}$. Proposition \ref{maillotprop} implies that
\begin{equation}
\label{star}
\wh{c}_t(\ov{S})\wh{c}_t(\ov{Q})=1+a(\wt{c}_t(\ov{\E}))=
1-p_a(t),
\end{equation}
where the subscript $t$ denotes the corresponding Chern polynomial.
Multiplying both sides of (\ref{star}) by $1+p_a(t)$ and using the
properties of multiplication in ${\cal A}$ gives the equivalent form
\begin{equation}
\label{star1}
\wh{c}_t(\ov{S})*\wh{c}_t(\ov{Q})*(1+p_a(t))=1.
\end{equation}

We now note that the {\em harmonic number generating function}
\[
\dis
\sum_{i=0}^{\infty}\H_it^i=\frac{t}{1-t}+\frac{1}{2}\frac{t^2}{1-t}+
\frac{1}{3}\frac{t^3}{1-t}+\cdots=
\frac{\log(1-t)}{t-1}.
\]
It follows that
\[
\dis
p_a(-t)=\sum_{i=0}^{\infty}\H_ip_i(y)t^{i+1}=
t\sum_{j=1}^s\sum_{i=0}^{\infty}\H_i(y_jt)^i=
-t\sum_{j=1}^s\frac{\log(1-y_jt)}{1-y_jt}
\]
and thus
\[
\dis
p_a(t)=t\sum_{j=1}^s\frac{\log(1+y_jt)}{1+y_jt}.
\]
Substituting this in equation (\ref{star1})
 gives relation ${\cal R}_2$.
\endproof

\medskip

 Theorem \ref{presentation} shows that the relations in the Arakelov
Chow ring of $G$ are the classical geometric ones perturbed by a
new `arithmetic factor' of $1+p_a(t)$. While this factor is
closely related to the power sums $p_i(\ov{Q})$,
 the most natural basis of
symmetric functions for doing calculations in $CH(G)$ is the basis
of Schur polynomials (corresponding to the Schubert classes; see for
example [F], \S 14.7). 
The arithmetic analogues of the special Schubert classes
involve the power sum perturbation above; multiplication formulas
are thus quite complicated (see [Ma]). 

 In geometry the Chern roots $x_i$ and $y_j$ all `live' on the complete
flag variety above $G$. There are certainly natural 
 line bundles on the flag variety whose first Chern classes
 correspond to the roots in Theorem \ref{presentation}. However on
 flag varieties the situation is more complicated and our knowledge
is not as complete. We refer the reader to [T] for more details.

\end{document}